# A Fractal Perspective on Scale in Geography

Bin Jiang and S. Anders Brandt

Faculty of Engineering and Sustainable Development, Division of GIScience
University of Gävle, SE-801 76 Gävle, Sweden
Email: bin.jiang@hig.se, sab@hig.se

*(Draft: September 2015, Revision: June 2016)*

**Abstract**
Scale is a fundamental concept that has attracted persistent attention in geography literature over the past several decades. However, it creates enormous confusion and frustration, particularly in the context of geographic information science, because of scale-related issues such as image resolution, and the modifiable areal unit problem (MAUP). This paper argues that the confusion and frustration arise from traditional Euclidean geometric thinking, with which locations, directions, and sizes are considered absolute, and it is now time to revise this conventional thinking. Hence, we review fractal geometry, together with its underlying way of thinking, and compare it to Euclidean geometry. Under the paradigm of Euclidean geometry, everything is measurable, no matter how big or small. However, most geographic features, due to their fractal nature, are essentially unmeasurable or their sizes depend on scale. For example, the length of a coastline, the area of a lake, and the slope of a topographic surface are all scale-dependent. Seen from the perspective of fractal geometry, many scale issues, such as the MAUP, are inevitable. They appear unsolvable, but can be dealt with. To effectively deal with scale-related issues, we present topological and scaling analyses illustrated by street-related concepts such as natural streets, street blocks, and natural cities. We further contend that one of the two spatial properties, spatial heterogeneity is de facto the fractal nature of geographic features, and it should be considered to the first effect among the two, because it is global and universal across all scales, which among the practitioners of geography should receive more attention.

## 1. Introduction
Scale is an important, fundamental concept in geography, yet it has multiple definitions or meanings, some of which seem to be contradictory. Among the various definitions (Lam 2004), *map scale* is the most commonly used, referring to the ratio of distance on a map to the corresponding distance on the ground. Scale is also closely related to map generalization for selectively representing things on the Earth's surface on a map and it can refer to the pixel size of an image, i.e. *resolution*. An image with small pixels has *high resolution*, while one with big pixels has *low resolution*. In this regard, scale is synonymous with the level of detail of an image, which is closely related to the notions of scaling up and down (Wu et al. 2006, Kim and Barros 2002) for translating statistical inference and reasoning from one scale to another. Scale is also commonly used to refer to the *scope or extent* of a study area. A large scale of study area (such as a country), if mapped, implies a small-scale map, whereas a small scale of study area (such as a city), if mapped, implies a large scale map. Obviously, confusion and frustration arise from multiple, seemingly contradictory meanings, and how to translate statistical inferences across scales. On the other hand, the confusion and frustration make scale even more interesting and challenging. In addition to the quantitatively defined scales, there are other qualitatively defined scales, such as micro-, meso- and macro-scales, and local, regional, and global scales.

The concept of scale has generated extensive literature over the past two decades (e.g., Sheppard and McMaster 2004), along with emerging geospatial technologies including geographic information science and remote sensing (e.g., Tate and Atkinson 2001, Weng 2014). Between 1997 and 2014, Goodchild and his colleagues produced eight publications with *scale* in the titles, including two books (Quattrochi and Goodchild 1997, Zhang et al. 2014). There are of course numerous other writings in the literature where scale and scale related issues such as the modifiable areal unit problem (MAUP) have been of persistent interest and challenge in geography and in geographic information science in particular. However, previous discussions are usually constrained to Euclidean geometry because all the meanings of scale in geography are about sizes in a ratio, or an absolute value related to geographic features or their representations. As a consequence, the major concern surrounding scale is how it affects geospatial data collection and analysis results with respect to accuracy and reliability. This is understandable because maps are initially produced for depicting and measuring things on the Earth's surface. Unfortunately, most geographic features are not measurable, or the measurement is scale-dependent because of their fractal nature (Goodchild and Mark 1987, Batty and Longley 1994, Frankhauser 1994, Chen 2011b). For example, the length of a coastline, the area of a lake, and the slope for a topographic surface are all scale-dependent, so they should not be considered absolute. Unfortunately, to a large extent our fundamental thinking on scale issues so far has been based on Euclidean geometry.

Scale in fractal geometry (Mandelbrot 1982), as well as in biology and physics (Bonner 2006, Jungers 1984, Bak 1996), is primarily defined in a manner in which a series of scales are related to each other in a scaling hierarchy. For example, a coastline is a set of recursively defined bends, forming the scaling hierarchy of far more small bends than large ones (Jiang et al. 2013). Therefore, a new definition of fractal could explicitly be based on the notion of far more small things than large ones (Jiang and Yin 2014, Jiang 2015a), in analogy with Christaller's (1933) central place theory (cf. Chen, 2011a) where there are many small villages but few large cities. Another definition of scale is simply the measuring scale, ranging from smallest to largest, to measure both Euclidean and fractal shapes. This measuring scale, from an individual rather than a series point of view, is equivalent to image resolution or map scale. This measuring scale makes many geographers believe fractal geometry can be a useful technique for dealing with scale issues. However, this view of fractal geometry is dubious. Fractal geometry is not just a technique but could also offer a new paradigm or new worldview that enables us to see surrounding things differently. Fractal geometry is a science of scale because it involves the universal scaling pattern across all scales from smallest to largest. On the contrary, geography dominated by traditional Euclidean geometry focuses on a few scales for measuring individual sizes.

This paper aims to advocate fractal thinking as a way to effectively deal with scale issues in geography and geospatial analysis. We think that mainstream views on scale of practitioners in geography, as briefly reviewed above, are not in line with the same concept in other sciences such as physics, biology and mathematics. In spite of the fact that fractal geometry has been intensively studied in geography (Goodchild and Mark 1987, Frankhauser 1994, Batty and Longley 1994, Chen 2011b), the fundamental way of thinking of most geographers while dealing with scale issues is still Euclidean. This situation still exists more than forty years after fractal geometry was established. On the other hand, this situation is understandable, because, with the development of geospatial technology, measurement with high accuracy and precision has been a major concern. To measure things, we need Euclidean geometry, whereas to develop new insights into structure and dynamics of geographic features, we need fractal geometry.

Section 2 introduces fractal geometry, in particular the underlying way of thinking, and put it in comparison with the Euclidean counterparts. Based on fractal geometry or fractal thinking, Section 3 presents several fallacies or scale related issues in geography, such as the conundrum of length and MAUP. To avoid these scale related issues, Section 4 illustrates street-based topological and scaling analyses that enable us to see the underlying scaling patterns. Finally Section 5 discusses two spatial properties that are closely related to the notion of scale, and summaries our major points towards the fractal perspective on scale.

## 2. Fractal geometry and the underlying way of thinking

Euclidean geometry is all about regular shapes, such as line segments, rectangles and circles, together with their locations, directions, and sizes, whereas fractal geometry deals with irregular shapes such as snowflakes and coastlines, with the underlying scaling property of far more small things than large ones. In this section, a line segment and the Koch curve are used as working examples to illustrate the two geometries, and the two different underlying ways of thinking.

### 2.1 Koch curve versus line segment

The Koch curve is one of the first fractals, invented by the Swedish mathematician Helge von Koch in 1904. If we have a line segment of one unit (or initiator), divide it into three thirds and replace the middle third by two sides of an equilateral triangle, the process results in a curve of four segments at iteration 1 (or generator). For each of the four segments, which is called generator, the same process of division and replacement is repeated, leading to a curve of 16 segments at iteration 2, and further recursively, a curve of 64 segments at iteration 3 (Figure 1). Theoretically, the same process can go on forever, leading to what is commonly known as the Koch curve, which is self-similar, with the similarity or scaling ratio of 1/3. What is interesting is that as the scale (or the line segment) gets smaller and smaller (1/3, 1/9, 1/27, …), the resulting curve becomes longer and longer (4/3, 16/9, 64/27, …), eventually having an infinite length. From the example of the Koch curve, the two meanings of scale in fractal geometry can be seen. The first is that the measuring scale (1/3, 1/9, 1/27, …) is a series. The second is that any curve containing far more short segments than long ones is also a series. It is essential to take a recursive perspective to clearly see the second meaning. Taking the curve of 64 segments as a whole, the black curve at iteration 3 in Figure 1 represents a non-recursive perspective, while the blue curve denotes a recursive version with 85 segments (i.e., 1 + 4 + 16 + 64). If the head/tail breaks classification method (Jiang 2013) is applied on the 85 segments, the head is the green curve (head 1 in Figure 1). By recursively applying head/tail breaks on head 1, we get head 2 and then 3. In other words, the orange curve contains one segment of one unit, and four segments of 1/3; the green curve contains one segment of one unit, four segments of 1/3, and 16 segments of 1/27; and the blue curve contains one segment of one unit, four segments of 1/3, 16 segments of 1/27, and 64 segments of 1/81.

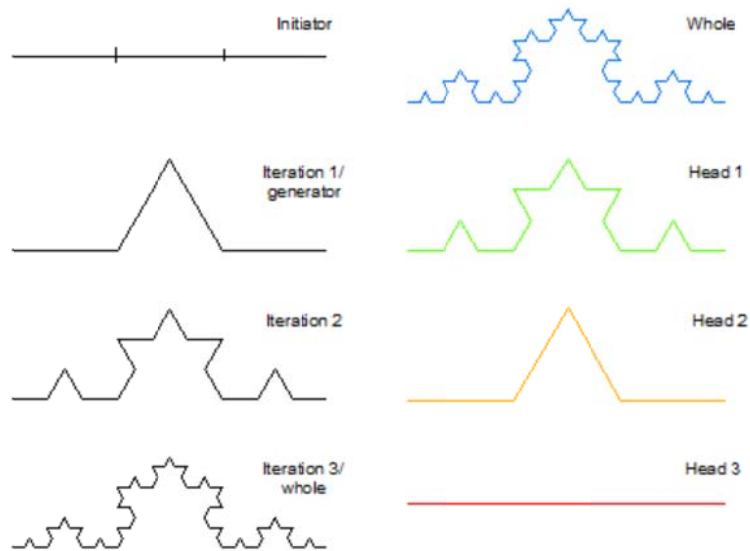

Figure 1: (Color online) Generation and decomposition of the Koch curve
(Note: The black curves are seen in a non-recursive view or Euclidean geometric view, while the color curves are in a recursive view or fractal geometric view. For example, the blue curve looks the same as the black curve at iteration 3, but the blue one is seen from the fractal geometric perspective, through which all the subsequent heads shown in color are embedded in the blue one. This recursive or fractal geometric perspective is based on head/tail breaks – a classification scheme and visualization tool for data with a heavy-tailed distribution (Jiang 2013, 2015a).)

In the course of generating the Koch curve, the initial segment of one unit is a simple regular shape. The Koch curve, on the other hand, looks complex, although the generation involves the simple repetition of division and replacement. First, the Koch curve is irregular and much more complex than the initial segment. Second, the Koch curve has an infinite length. Under the Euclidean geometry framework, anything is measurable, no matter how big or small. Why the Koch curve has an infinite length puzzled mathematicians for more than 100 years, until Mandelbrot (1967) solved the mystery. Geographic features such as coastlines bear the same property of the Koch curve. The length of a coastline is unmeasurable, or specifically, it is scale-dependent. In this way, the Koch curve and a coastline are essentially the same in terms of scale dependence. However, a coastline belongs to a statistical fractal with a limited scaling range, while the Koch curve is a strict fractal with an infinite scaling range. Therefore, fractal geometry offers a new worldview for viewing surrounding things such as trees, coastlines, and mountains.

**2.2 Fractal and Euclidean thinking**
Besides that Euclidean geometry considers regular simple shapes, and fractal geometry irregular complex shapes, there are more profound facets in how the two geometries differ (Mandelbrot and Hudson 2004) (Table 1). Euclidean geometry focuses on pieces or parts, while fractal geometry focuses on the whole. Euclidean geometry looks at individuals, while fractal geometry looks at patterns. This holistic or pattern view of fractal geometry implies a recursive view of seeing surrounding things. The Koch curve at iteration 3 (Figure 1) is just a Euclidean shape that consists of 64 segments of all the same scale of 1/27. Seen from the recursive or fractal geometric perspective, it becomes a fractal shape, involving 85 segments of four different scales with far more short scales than long ones.

Table 1: Comparison of Euclidean and fractal thinking

| Euclidean thinking | Fractal thinking |
| --- | --- |
| Regular shapes | Irregular shapes |
| Simple | Complex |
| Individuals | Pattern |
| Parts | Whole |
| Non-recursive | Recursive |
| Measurement ( = scale) | Scaling ( = scale free) |

Hence, fractal geometry offers a way of seeing our surrounding geography differently. Euclidean geometry mainly measures shapes (Euclidean shapes), directions, and sizes. Fractal geometry aims to see underlying scaling. Simply put, Euclidean geometry is used for one particular scale or a few scales, while fractal geometry aims for scale-free or scaling that involves all scales. The term *scale-free* is synonymous with scaling, literally meaning no characteristic mean for all sizes. This difference is very much like that between Gaussian and Paretian thinking (Jiang 2015b), which refer to more or less similar things (with a characteristic mean), and far more small things than large ones (without a characteristic mean, or scale-free), respectively. For example, a tree is better characterized by all sizes of its branches, or how the branches (scales) form a scaling hierarchy of far more small branches than large ones, rather than only by its height. It is fair to say that both Euclidean and fractal geometries aim to characterize things, but with different means; the former through measurement (at one scale), and the latter through scaling (across all scales). However, without individual Euclidean shapes, there would be no fractal pattern. It is scale that bridges individual Euclidean shapes and a fractal pattern. Without scale, there would be no fractal geometry. Scale plays the same important role in geography as in fractal geometry as a science of scale.

Fractal geometry is not just limited to patterns. It can also be applied to a set. For example, the set of numbers, 1, 1/2, 1/3, 1/4, …, 1/1000, constitutes a fractal because there are far more small numbers than large ones within the set, based on the definition of fractal using head/tail breaks classification method (Jiang and Yin 2014, Jiang 2015a). The 1,000 numbers are created by following Zipf's Law (Zipf 1949); the first number is 1, the second is just 1/2, the third number is 1/3, and so on. Hence, to a great extent a fractal pattern is more a statistical pattern that meets the scaling law than a visual pattern. Therefore, at the fundamental base, fractal thinking is not much different from Paretian thinking (McKelvey and Andriani 2005, Jiang 2015b) because both are concerned with the scaling pattern of far more small things than large ones.

### 3. Fallacies of scale in geography

Geographic features, both natural and manmade, are fractal in nature (Goodchild and Mark 1987, Batty and Longley 1994, Chen 2011b, Jiang 2015a), although within a limited scaling range. However, for the purpose of measurement, we have to assume geographic features as Euclidean or non-fractal. Because of this assumption, geographic features are considered measurable, or the measurement to be absolute. This is a common fallacy in geography. In addition, it is commonly considered that the slope for a topographic surface is computable. However, a tangent cannot be defined at any point of a coastline. The same is true for a topographic surface devoid of tangent planes (Mandelbrot 1982). The lack of a tangent line or plane implies that curve lengths and surface slopes are not measurable, yet people often treat these measurements as something absolute. In this section, we present some fallacies commonly seen in geography on the sizes of geographic features and on translating statistical results across scales.

#### 3.1 Coast length, island area, and surface slope

The length of a linear geographic feature such as a coastline is not measurable. To be more precise, the measurement depends on map scale or image resolution. As map scale exponentially increases, the coastline length tends to increase exponentially at the same pace (Richardson 1961). This conundrum of length, also commonly known as *paradox of length* (Steinhaus 1983) puzzled scientists for a very long time. Geographers attempted to solve this problem (Perkal 1966, Nystuen 1966) in order to measure geographic features on maps. However, the problem is essentially unsolvable because of the fractal nature of geographic features. Mandelbrot (1982) eventually uncovered the conundrum of length and further developed fractal geometry.

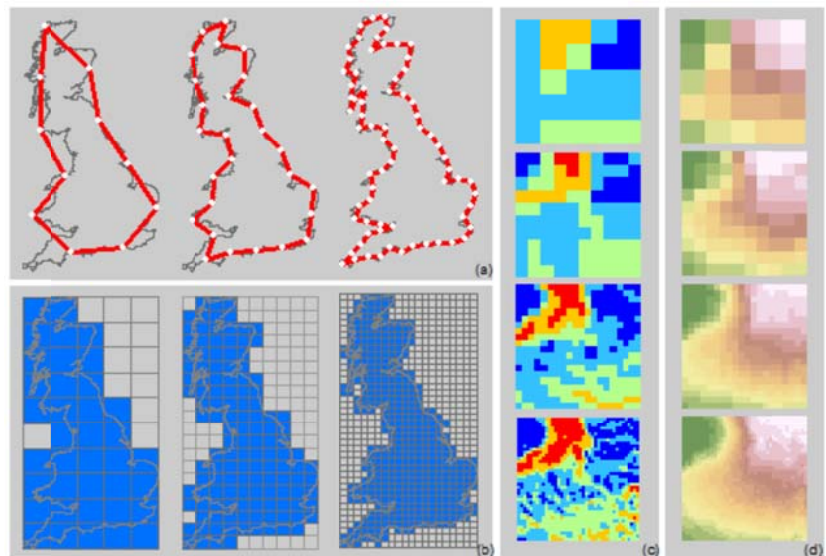

Figure 2: (Color online) Illustration of how measurement changes with scales
(Note: The length of the coastline changes with yardsticks (Panel a), the area of the island and the slope of the surface change with boxes and digital elevation model (DEM) resolutions (Panels b and c). The five slope classes between blue and red are respectively, 0–6.8, 6.8–12.4, 12.4–18.7, 18.7–26.9, and 26.9–40.9 degrees. Derived for the same area but with different resolutions of DEMs (Panel d), the slopes look very different with different resolutions, such as more areas of high slopes in higher-resolution DEMs.)

Although the conundrum of length is unsolvable, it can be effectively dealt with. It can also be expressed as such: As measuring scales decrease exponentially, the length of a coastline increases exponentially; see Panel (a) in Figure 2, where the length of the British coastline increases as the measuring yardstick decreases. The change of scales (either map scale or measuring scale) on the one hand and the change of the length on the other meet a linear relationship at logarithm scales, which is commonly shown in Richardson's plot (Richardson 1961). This simple statistical relationship enables us to predict the length of a coastline at different scales of maps, or to transform the coastline length across different scales of maps. However, such a transformation is only possible within the scaling range in which the simple linear relationship holds.

The kind of scale problem with a coastline also applies to areal geographic features such as islands and lakes. The area of an island depends on map scale or image resolution. In Panel (b) (Figure 2), a raster cell is marked and its area included as soon as any land (or country border line) is within its boundary. Hence, as the image resolution increases, the size of the UK decreases and then asymptotically converges to the true area value. Because of its fractal nature, the size of an island is transformable from one scale to another within the scaling range. The scale-dependent lengths and areas can have a similar effect on some second-order measurements, such as density defined on areas and kernel-density estimation based on a set of distance parameters.

Panel (c) in Figure 2 further illustrates that surface slope depends on the resolution of DEMs. In theory, the calculation of slope is based on the assumption that the topographic surface is differentiable, or with tangent planes. However, this assumption is not true because topographic surfaces are fractal. In the figure, the different slope maps are shown, representing DEMs based on light detection and ranging from different flight heights. While the elevation range for the four DEMs is more or less identical (Panel d in Figure 2), the derived slope range is considerably different between the DEMs. The highest resolution DEM has a range of 0 – 40.9 degrees, while the lowest resolution DEM has a range of 5.3 – 22.6 degrees. This infers that: first, slope values are scale dependent; second, with decreasing resolution, the slope values tend to converge to one single value. This behavior can thus be directly linked to the fractal nature of terrain surfaces. To further illustrate this, a series of slope maps (ranging from 1 to 4 m resolution) for an area are shown as a 3D-surface, where the total area of each slope class is plotted against the 1-degree wide slope classes (Figure 3). For the finest resolution (or scale), the slope-value distribution is the most heterogeneous (i.e. a flatter and more evenly distributed area-slope curve), whereas if only the coarsest resolution (or scale) is considered, the slope-value distribution become less heterogeneous.

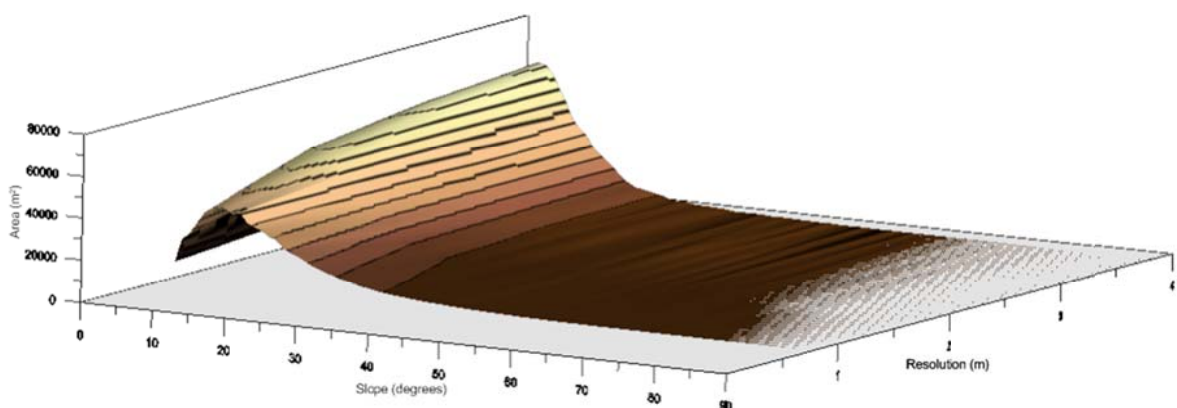

Figure 3: (Color online) Illustration of how scale change can be related to change of slopes

### 3.2 The MAUP
Geographic spaces are usually organized hierarchically in a nested manner in which small units are contained in large ones, i.e. far more small units than large ones. For example, the 50 US states constitute the first administrative level, followed by thousands of counties, and further by many more cities and towns. A problem arises when putting measurements into different areal units for statistical analysis, i.e., the analysis results vary from one level of unit to another. This is called the MAUP. This MAUP phenomenon was first discovered by Gehlke and Biehl (1934) when they studied crime rate and average income, and found that the correlation coefficient increases when the data was aggregated into large areal units. The term MAUP was later coined by Openshaw and Taylor (1979), who studied the issue systematically, and since then the MAUP has created a large body of literature and appears in almost every book on spatial analysis (e.g., Wong 2004, Fotheringham and Rogerson 2009), including books on scales (Wu et al. 2006, Lock and Molyneaux 2006). These books demonstrate how the MAUP has drawn the attention of especially archaeologists and landscape ecologists concerned with how to translate statistical results across different scales, namely scaling up and down. This is because the crux of the MAUP lies in the translation of statistical inference from an aggregated level to an individual level, or from large to small aggregated levels, commonly known as an ecological fallacy (Robinson 1950, King 1997).

In theory, statistical inference and reasoning cannot be simply translated from one level to another or from an aggregated scale to an individual one. For the sake of simplicity, Figure 4 illustrates both scale and zoning effects of the MAUP, using an example of simple statistic of the unemployment rate. The rate varies from one level to another (the scale effect), and from one configuration to another (the zoning effect). The most accurate indicator is with the smallest unit. A typical example of the zoning effect is gerrymandering, which establishes an advantage for a particular political party by manipulating districts to create advantageous voting areas. Gerrymandering could dramatically reverse initial election results due to the zoning effect. It is certainly of general scientific interest to learn how analysis results obtained at one scale would apply to other scales of geospatial data or other extents (larger or smaller) of study areas. However, such a translation across scales is deemed to be flawed, especially across scales of a large gap in the hierarchy by either scaling up or down.

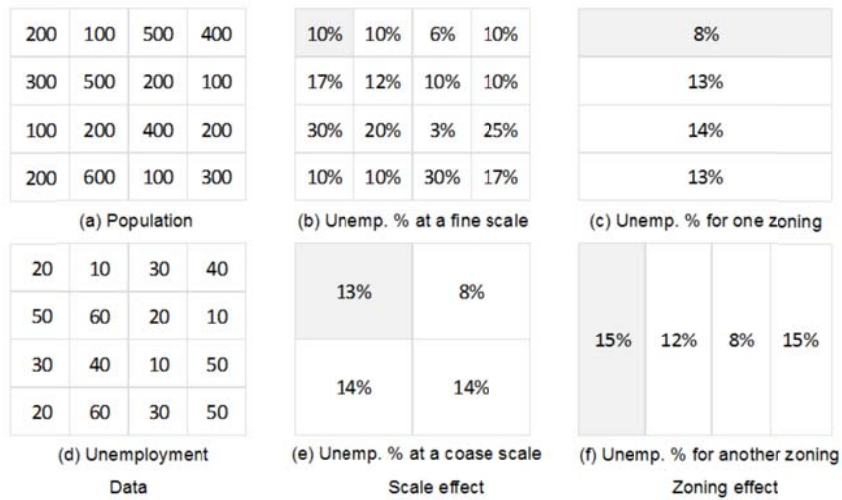

Figure 4: Illustration of the MAUP using unemployment rate as an example
(Note: The original data are shown in panels (a, d), while the derived percentages in panels (b, c, e, f). For example, the four percentages in the shaded cells are as follows: 10% = 20/200, 8% = (20 + 10 + 30 + 40)/(200 + 100 + 500 + 400), 13% = (20 + 10 + 50 + 60)/(200 + 100 + 300 + 500), and 15% = (20 + 50 + 30 + 20)/(200 + 300 + 100 + 200).)

The MAUP also applies to other kinds of data. Take the example of a large number of students who are put into different groups according to performance in different courses, where the group behavior cannot directly be translated into individuals. It might be possible only by sacrificing certain reliability and if the groups are very homogeneous. If the units are small enough, and individuals in each group are likely to be homogeneous, translation of statistical results from groups to individuals might be possible in practice. However, such a translation across scales is hardly possible if the groups are very heterogeneous. For example, predicting weather conditions between two consecutive weeks is much harder than between two consecutive days, which is even harder than between two consecutive hours.

A common problem is that in some research areas such as demography and ecology, there is usually no individual-based data, and all data are in aggregated formats. Statistical results apply to the particular scale in which data is collected. Translation of the results into other scales, either scaling up or down, is impossible unless there is a relationship between statistical inferences and scale changes. For example, for a curve length, there exists such a relationship, i.e., the ratio of the change in length to the change in scale at logarithmic scales is constant. If there is no such a relationship, keep the inference only to that particular scale. The advent of big data has changed the situation, in which data is individually based and units can be objectively derived. One example of this is how cities can be more objectively defined (see Section 4). Therefore, a good solution to effectively dealing with the MAUP is to conduct frame-independent spatial analysis (Tobler 1989). Therefore, the next section presents topological and scaling analyses to effectively deal with scale issues.

## 4. Topological and scaling analyses to deal with scale issues

Euclidean geometric thinking aimed for measuring things would inevitably result in errors, uncertainties, or scale effects in general. To effectively deal with scale issues, we must shift our thinking from geometric details to topological relationships, and from measurement to scaling analysis. Topological thinking concentrates on relationships of things and enables us to see the underlying scaling pattern of far more less-connected things than well-connected ones, in analogy with far more small branches than big branches in a tree or a river network. We exemplify this through some street-related concepts to elaborate on topological and scaling analyses that are free of scale effects.

### 4.1 Natural streets

Euclidean geometry has critical limitations in developing analytical insights into geospatial data, but this viewpoint has not always been well received. For example, a street network is a graph of individual street segments or junctions. This graph representation of street networks is Euclidean-oriented in terms of junction locations and distances of street segments. Despite its usefulness in computing routes and distances, the geometric representation offers little analytical insights into the underlying structure. There are similar segments in terms of length or similar junctions in terms of degree of connectivity; in other words, a boring structure. Therefore, let us shift from the geometrical details of segments and junctions to the topological structure of individual streets.

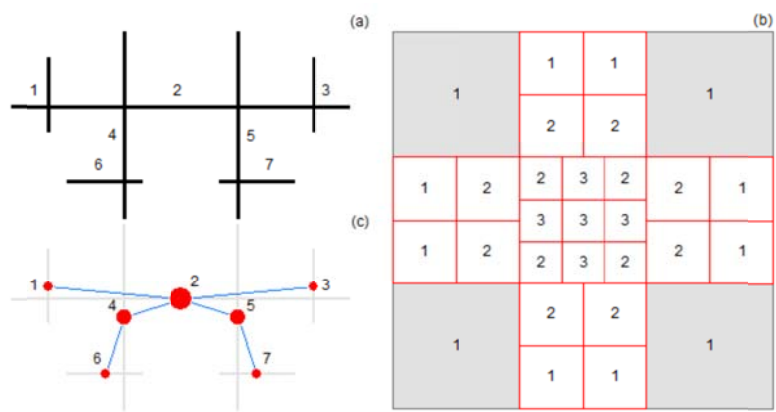

Figure 5: (Color online) Illustration of concepts of natural streets and a natural city
(Note: The seven natural streets, indicated by numbers in Panel (a), are created from 19 street segments, while the one natural city is created from the 25 city blocks, which are red in Panel (b). There are far more less-connected streets than well-connected ones, and far more small street blocks (or city blocks) than large ones. In Panel (c), the seven streets and their six intersections are converted into the seven nodes and six links of the connectivity graph; the dot sizes show the degrees of connectivity. The 25 city blocks and the four large street blocks (or field blocks in gray) in Panel (b) constitute the fictitious country as a whole. The numbers shown in the blocks indicate how topologically far the blocks are from the country's border, or the border number.)

So called 'natural streets' are created from adjacent street segments with good continuity or with the least deflection angles or with the same names (Jiang et al. 2008, Jiang and Claramunt 2004). These streets are naturally defined and can be a basic unit for spatial analysis through a graph representation, or a connectivity graph, in which the nodes and links represent individual natural streets and their intersections (see Panels (a, c) in Figure 5). This topological representation provides an interesting structure, involving all kinds of streets in terms of both length and degree of connectivity. In other words, interconnected streets behave as fractal when seen from the topological perspective, as there are far more less-connected streets than well-connected ones. Thus, natural streets can be a more meaningful unit than arbitrarily imposed areal units. Natural streets, or their topological representations, suffer less from scale effect because geometric details such as accuracy and precision play a less important role. What matters is relationship. The topological view enables us to see the underlying scaling pattern, where we can assign point-based data into individual natural streets, rather than to any modifiable areal units for spatial analysis.

### 4.2 Street blocks and natural cities

Street blocks also demonstrate the scaling property of far more small things than large ones. The street blocks refer to the minimum rings or cycles, each of which consists of a set of adjacent street segments. Obviously, a country's street network usually comprises a large amount of street blocks (Jiang and Liu 2012). The street blocks are the smallest unit and are defined from the bottom up, rather than imposed from the top down by authorities. The street blocks are smaller than any administratively or legally imposed geographic units. They can be automatically extracted from all kinds of streets, including pedestrian and cycling paths. It is understandable that the street or city blocks defined by authorities are just a subset of the automatically extracted street blocks. For example, the number of London census output areas, which is the smallest census unit in the UK, is just half that of the street blocks that can be extracted from the OpenStreetMap databases.

Topological analysis of the street blocks begins with defining the border number, which is the topological distance far from the outermost border of a country. The border is not a real country border, but consists of the outermost street segments of the street network. Those blocks adjacent to the border have border number one, and those adjacent to the blocks of border number one have border number two, and so on (see Panel (b) in Figure 5). All the blocks are assigned a border number, indicating how far they are from the outermost border. Interestingly, the block(s) with the highest border number constitutes the topological center of the country (Panel (b) in Figure 5). In the same way, we can take all city blocks as a whole to define the topological center as the city center. The topological center differs from the geometric center, or the central business district that is commonly the city center.

The scaling property of far more small blocks than large ones enables us to define the notion of natural cities emerging from a large amount of heterogeneous street blocks. All the street blocks are inter-related to form a whole. The whole can be broken into the head for those above the mean, and the tail for those below the mean, as shown in the head/tail breaks classification method described in Jiang (2013). Those small street blocks in the tail constitute individual patches called natural cities. See Panel (b) in Figure 5 for an example of a natural city. Natural cities are defined from the bottom up. A large amount of street blocks collectively decides a mean value as a cut off for the city border. The head/tail breaks can recursively continue to derive patches within individual natural cities. In other words, all city blocks within a natural city are considered a whole, and those below the mean value in the tail (high-density clusters) are considered hotspots of the natural city. It is essentially the fractal nature of street structure, or the scaling property of far more small blocks than large ones, that make the natural cities definable. We can assign point-based data into city blocks, or natural cities, rather than any modifiable areal units for spatial analysis.

In summary, as many modifiable areal units are imposed by authorities or images from the top down, such as administrate boundaries, census units, and image pixels, it is inevitable that these units or boundaries are somehow subjective. They were defined mainly during the small-data era for the purpose of administration and management, but are still used in the big-data era. However, objectively defined units such as natural streets, street blocks, and natural cities should therefore be better alternatives for scientific purposes, and they reflect the new ways of thinking about data analytics in this big data era.

## 5. Discussion and summary

The concept of scale is closely related to spatial heterogeneity, one of the two fundamental spatial properties. The other property is spatial dependence, or auto-correlation, which has been formulated as the first law of geography: *Everything is related to everything else, but near things are more related than distant things* (Tobler 1970). The Tobler's law implies that near and related things are more or less similar. Therefore, spatial variation for near and related things is mild rather than wild in terms of Mandelbrot and Hudson (2004). However, spatial heterogeneity is about far more small things than large ones, or with wild rather than mild variation. The two spatial properties are closely related and can be rephrased as such. There are far more small things than large ones in geographic space – spatial heterogeneity, but near and related things are more or less similar – spatial dependence. In this regard, spatial heterogeneity appears to be the first-order effect being global, while spatial dependence is the second-order effect being local. Therefore spatial heterogeneity provides a larger picture across all scales ranging from smallest to largest, while spatial dependence is a more local pattern.

As a fundamental concept in geography, scale has been a major concern for geospatial data collection and analysis, not only in geography and geographic information science, but also in ecology and archaeology. Although Goodchild and Mark (1987, p. 265) concluded that *"fractals should be regarded as a significant change in conventional ways of thinking about spatial forms and as providing new and important norms and standards of spatial phenomena rather than empirically verifiable models"*, this paper suggests fractal geometry thinking finally should become a new paradigm, rather than a technique recognized in the current geographic literature. To effectively tackle scale issues in geospatial analysis for better understanding geographic forms and processes, this fractal or recursive perspective is essential. Many people tend to think that a cartographic curve is just a collection of line segments – Euclidean way of thinking, but actually it consists of far more small bends than large ones – fractal way of thinking. After comparing various definitions of scale in geography and cartography, as well as in fractal geometry, we note that the definitions of scale in geography are very much constrained by Euclidean geometry, for measuring geographic features rather than for illustrating the underlying scaling pattern. Most geographic features are inevitably fractal, so their sizes are unmeasurable or scale-dependent. We must be aware of scale effects in measuring geographic features, and that their sizes change as the measuring scale changes. The measurement is a relative indicator, rather than something absolute.

Not only their sizes, but also statistical inferences on geographic features are scale dependent. Statistical reasoning cannot be translated across scales, or from an aggregate scale to individual ones. Given these circumstances, we must adopt fractal thinking for geospatial analysis involving all scales rather than a single scale or a few scales. We must examine if there are far more small things than large ones, rather than measuring individual sizes. We must also determine if, or how, the locations are related, rather than measuring absolute locations. However, unlike their sizes, statistical inferences on geographic features have not been found to hold a simple relationship with scales within a scaling range. This certainly warrants further research in the future.

**Acknowledgment**
XXXXXX